\documentclass[11pt]{article}
 \hoffset= - 0.75 in

\usepackage{amsfonts, amsmath, amssymb, amsthm}
\usepackage[xdvi]{graphics}
\usepackage{latexsym}
\newtheorem{thm}{Theorem}[section]
\newtheorem{defn}[thm]{Definition}

\newtheorem{conj}[thm]{Conjecture}

\newtheorem{lem}[thm]{Lemma}
\newtheorem{cor}[thm]{Corollary}
\newtheorem{prop}[thm]{Proposition}

\newtheorem{prob}[thm]{Problem}

\newtheorem{Example}[thm]{Example}
\newenvironment{exmp}{\begin{Example}\rm}
                {\mbox{}~\hfill~\end{Example}}

\newcommand{\ZZ}{\mathbb{Z}}
\newcommand{\NN}{\mathbb{N}}
\newcommand{\RR}{\mathbb{R}}

\begin{document}

\title{Markov Bases of Binary Graph Models}
\author{Mike Develin \\ American Institute of Mathematics\\ 360 Portage Ave.
\\ Palo Alto, CA 94306-2244, USA \\ develin@post.harvard.edu
\and Seth Sullivant \\ Department of Mathematics \\ UC-Berkeley \\
Berkeley, CA 94720-3840, USA \\ seths@math.berkeley.edu}

\maketitle

\begin{abstract}
This paper is concerned with the topological invariant of a graph
given by the maximum degree of a Markov basis element for the
corresponding graph model for binary contingency tables.  We
describe a degree four Markov basis for the model when the
underlying graph is a cycle and generalize this result to the
complete bipartite graph $K_{2,n}$. We also give a combinatorial
classification of degree two and three Markov basis moves as well
as a Buchberger-free algorithm to compute moves of arbitrary given
degree.  Finally, we compute the algebraic degree of the model
when the underlying graph is a forest.

\smallskip

\end{abstract}

\section{Introduction and Definitions}

The study of multidimensional tables and their marginals is of
central importance whenever one wishes to make inferences based on
statistical samples.  In general, one is presented with a
nonnegative integral table of data of size $d_1 \times \cdots
\times d_n$ and a simplicial complex $\Delta$ on $\{1, \ldots,
n\}$ which encodes the specific marginals we would like to
compute; this is called a hierarchical model.

 Certainly, the oldest example of such a model
is the case of computing the row and column sums of a matrix.  If
the matrix is a square $m \times m$ matrix and we require that all
the row and column sums have the same value $k$, one calls such a
matrix a semi-magic square with magic sum $k$. Here, our
simplicial complex consists of two isolated points.

In the more general statistical situations, each node of the
simplicial complex corresponds to a feature of a population sample
(e.g. eye color) and the levels of the table correspond to
different states of the feature (e.g. green, brown, hazel, etc.).
The faces of the simplicial complex are intended to model
interactions between the features.  One of the most fundamental
questions in statistical analysis is:  do the data appear to be
satisfied by a given model?  One way to test the hypothesis is to
compare the sample data to the maximum likelihood estimate using
the $\chi^2$ or $G^2$ statistic.  One problem with applying this
approach directly is that the data is always an integral table,
while the maximum likelihood estimate almost never is.  As a
result, there might be no integral table with the same marginals
which has a small $\chi^2$ statistic.  This problem is especially
dramatic when analyzing the large, sparse data sets which occur in
real world situations (e.g. census data). To remedy this
situation, one may attempt to decide whether or not a model fits
the data by comparing statistics of the data table with statistics
of random integral tables with the same marginals. If the
statistic of the table of data is exceptional one could hope to
conclude that the table was also exceptional. For example, if the
$\chi^2$ statistic of the table of data was exceptionally small,
one could conclude that the data did, in fact, fit the model.

We are now left with the problem of generating random integral
tables from the set of all nonnegative integral tables with fixed
marginals.  One solution is to perform a random walk over the set
of all nonnegative integral tables with given fixed marginal. Such
a random walk can be taken by first finding a suitable set of
``moves'' (these are tables with integral entries which have all
their marginals equal to zero) and randomly adding moves to some
starting table.  It is at this point in the story that
computational commutative algebra enters the picture: finding such
a set of moves is equivalent to finding a generating set for the
associated toric ideal.  For a detailed introduction to the
connections between toric algebra and multidimensional contingency
tables, see \cite{DS}, \cite{HS}, and \cite{St2}.

In this paper, we restrict attention to tables where $d_i = 2$ for
all $i$ and for which the underlying simplicial complex is a
graph; that is, we compute only two- and one-way marginals of our
binary table. We refer to such models as binary graph models.
These are generally not the usual graphical models studied so
frequently in statistics, where the simplicial complex consists of
the cliques of the underlying graph $G$ \cite{L}.  Our notion of
graph model coincides with the more familiar graphical model if
and only if the graph has no three-cycle.  Now, we will give two
formal presentations of the objects of interest in this paper.

Let $G$ be a graph on the $n$-element vertex set $[n] := \{1,2,
\ldots,n\}$ with edge set $E(G)$. Denote by $Iso(G)$ the set of
isolated vertices of $G$.  For each edge $\{j,k\}$ of the graph
$G$ consider the linear transformations $\pi_{j,k}$

$$ \pi_{j,k} : \ZZ^{2^n} \longrightarrow \ZZ^{4}$$
$$  e_{i_1, \ldots i_n} \longmapsto  e_{i_j,i_k}$$

\noindent and for each isolated vertex $k$ of $G$ consider the
linear transformations $\pi_k$

$$ \pi_k : \ZZ^{2^n} \longrightarrow \ZZ^2 $$
$$ e_{i_1 \ldots i_n} \longmapsto  e_{i_k}.$$

\noindent  We think of the maps $\pi_{j,k}$ as computing the 2-way
marginal of a $2 \times \cdots \times 2$ table corresponding to
the edge $\{j,k\}$ and the maps $\pi_{k}$ as computing the 1-way
marginal of a $2 \times \cdots \times 2$ table corresponding to
the vertex $k$.  We define the map $\pi_G$ by taking all the
marginal computations induced by a given graph as

$$\pi_G:  \ZZ^{2^n} \longrightarrow \bigoplus_{\{j,k\} \in E(G)}
\ZZ^{4} \bigoplus_{k \in Iso(G)} \ZZ^{2}$$

$$ v \longmapsto \bigoplus_{\{j,k\} \in E(G)}
\pi_{j,k}(v)  \bigoplus_{k \in Iso(G)} \pi_k(v).
$$

We say that $\pi_G$ is the map which computes the marginals of a
$2 \times \cdots \times 2$ table according to the graph $G$.  The
matrix which represents this linear transformation will be denoted
$A_G$ and the polytope which is the convex hull of the columns of
$A_G$ is denoted $P_G$ where we consider the columns of $A_G$ as
vectors in $\RR^d$ for an appropriate $d$.  A move for $G$ is an
element of the integral kernel of $\pi_G$; that is, a move is an
integral table which does not change the $G$-marginals of a table
it is added to.  In general we are interested in sets of moves
with special properties.

\begin{defn}
A finite subset of moves $B \subset \ker_\ZZ( \pi_G)$ is called a
{\bf Markov basis} for the graph $G$ if for every pair of
nonnegative integral tables $v_1, v_2 \in \NN^{2^n}$ with the same
$G$-marginals $\pi_G(v_1) = \pi_G(v_2)$, there is a sequence of
moves $\{u_i\}_{i=1}^l \subset \pm B$ such that

$$
v_1 + \sum_{i=1}^l u_i= v_2
$$
and
$$ v_1 + \sum_{i =1}^j u_i  \in \NN^{2^n} \quad \mbox{ for all } 1
\leq j \leq l .$$
\end{defn}

There is also a much shorter and more algebraic way to arrive at
this definition.  Recall that $G$ is a graph on the $n$-element
vertex set $[n] := \{1,2, \ldots,n\}$ with edge set $E(G)$, and
isolated vertices $Iso(G)$. Consider the map of polynomial rings

$$\phi_G:  \mathbb{C}[p_{i_1,\ldots i_n} | i_j \in \{0,1\}] \longrightarrow
\mathbb{C}[t^{(j,k)}_{i_j, i_k}, t^{(l)}_{i_l} | \{j,k\} \in E(G)
\mbox{ and } l \in Iso(G)]
$$

$$ p_{i_1 \ldots i_n} \longmapsto \prod_{\{j,k\} \in E(G)}
t^{(j,k)}_{i_j,i_k} \prod_{l \in Iso(G)} t^{(l)}_{i_l}.$$

\noindent  The object of interest in this paper is the ideal
denoted $I_G = \ker(\phi_G)$, which we call the ideal of a binary
graph model.  It is a {\bf toric ideal}: a prime ideal generated
by monomial differences whose leading and trailing terms have
disjoint support.

Markov bases and toric ideals are connected by the following
fundamental theorem.

\begin{thm}\cite{DS} \label{thm:mbgs}
A finite subset of moves $B = \{u_i\}_{i=1}^l \subset
\ker_\ZZ(\pi_G)$ is a Markov basis for $G$ if and only if the set
of binomials $\{p^{u_i^+} - p^{u_i^-} \}$ is a generating set for
$I_G$.
\end{thm}

Here, we write $u_i = u_i^+ - u_i^i$ as the difference of two
positive vectors of disjoint support.  In light of Theorem
\ref{thm:mbgs}, we will use the expressions ``Markov basis for
$G$'' and ``generating set for $I_G$'' interchangeably throughout
this paper.  Similarly, we can interchange the words ``move'' and
``binomial'' whenever we are discussing the Markov
bases/generating sets of $I_G$.  These definitions are best
illustrated by a simple example.

\begin{exmp}
Consider the graph $G$ on four nodes with two edges $\{1,2\}$ and
$\{2,3\}$ and one isolated vertex $4$.  The map $\phi_G$ is a map
from a polynomial ring in sixteen variables to a polynomial ring
in ten variables.  It is the map of rings

$$\phi_G : \mathbb{C} [p_{ijkl} | i,j,k,l \in \{0,1\}]
\longrightarrow \mathbb{C} [r_{ij}, s_{jk}, t_l]$$
$$ p_{ijkl} \longmapsto r_{ij} \cdot s_{jk} \cdot t_l. $$

\noindent  On the other hand, the marginal map $\pi_G$ is given by
a $10 \times 16$ matrix $A_G$.  It is the matrix

$$A_G = \left( \begin{array}{cccccccccccccccc}
1 & 1 & 1 & 1 & 0 & 0 & 0 & 0 & 0 & 0 & 0 & 0 & 0 & 0 & 0 & 0 \\
0 & 0 & 0 & 0 & 1 & 1 & 1 & 1 & 0 & 0 & 0 & 0 & 0 & 0 & 0 & 0 \\
0 & 0 & 0 & 0 & 0 & 0 & 0 & 0 & 1 & 1 & 1 & 1 & 0 & 0 & 0 & 0 \\
0 & 0 & 0 & 0 & 0 & 0 & 0 & 0 & 0 & 0 & 0 & 0 & 1 & 1 & 1 & 1 \\
1 & 1 & 0 & 0 & 0 & 0 & 0 & 0 & 1 & 1 & 0 & 0 & 0 & 0 & 0 & 0 \\
0 & 0 & 1 & 1 & 0 & 0 & 0 & 0 & 0 & 0 & 1 & 1 & 0 & 0 & 0 & 0 \\
0 & 0 & 0 & 0 & 1 & 1 & 0 & 0 & 0 & 0 & 0 & 0 & 1 & 1 & 0 & 0 \\
0 & 0 & 0 & 0 & 0 & 0 & 1 & 1 & 0 & 0 & 0 & 0 & 0 & 0 & 1 & 1 \\
1 & 0 & 1 & 0 & 1 & 0 & 1 & 0 & 1 & 0 & 1 & 0 & 1 & 0 & 1 & 0 \\
0 & 1 & 0 & 1 & 0 & 1 & 0 & 1 & 0 & 1 & 0 & 1 & 0 & 1 & 0 & 1
\end{array} \right). $$

The polytope $P_G$ in $\RR^{10}$ has dimension 6 and has 10
facets.  These facets are indexed naturally by the rows of $A_G$
and the facet defining inequalities are given by $y_i \geq 0$ with
one inequality for each row.  The ideal $I_G$ has a Markov basis
consisting of  quadratic moves.  These are

$$ p_{0j0l_1}p_{1j1l_2} - p_{0j1l_1}p_{1j0l_2} \mbox{  with  } j,
l_1, l_2 \in \{0,1\}$$

\noindent and

$$ p_{i_1 j_1 k_1 l_1} p_{i_2 j_2 k_2 l_2} - p_{i_1 j_1 k_1 l_2} p_{i_2 j_2 k_2
l_1} \mbox{  with  } i_1, j_1, k_1, l_1, i_2, j_2, k_2, l_2 \in
\{0,1\}.$$

\noindent  These generators are also a Gr\"obner basis with
respect to the reverse lexicographic term order with $p_{0\cdots
0} \prec \cdots \prec p_{1 \cdots 1}$.
\end{exmp}

In general, when the underlying graph is a forest, the toric ideal
$I_G$ is relatively well understood.

\begin{thm}\label{thm:sqfr2}\cite{D, GMS, T}
The ideal $I_G$ is minimally generated by quadrics if and only if
$G$ is a forest.  In this case, the set of quadratic squarefree
binomials in $I_G$ forms a Gr\"obner basis with respect to the
reverse lexicographic term order with $p_{0 \cdots 0} \prec p_{0
\cdots 0 1} \prec \ldots \prec p_{1 \cdots 1}$.
\end{thm}

In this paper, we are primarily concerned with investigating
graphs which contain cycles. One fundamental question is to compute the
following invariant of a graph.

\begin{defn}
Let $G$ be a graph. The {\bf Markov width} $\mu(G)$ is the degree of the
largest minimal generator of the toric ideal $I_G$.
\end{defn}

Studying the Markov width of a graph is of fundamental importance
for statistical applications because it relates the complexity of
analyzing data to the complexity of the underlying graphical
structure.  Note that, Theorem~\ref{thm:sqfr2} states that the
graphs with $\mu(G)=2$ are precisely forests, which are certainly
topologically simple. Indeed, the Markov width of a graph $G$ is
topological in nature, by which we mean that $\mu(G)$ can only
decrease under the operations of vertex deletion and edge
contraction; we will show this in Section 4.  Since these
operations interact nicely with the toric ideals of the initial
and final graph, we will use the following definition throughout.

\begin{defn}
Let $G$ be a graph. By a {\bf minor} of $G$, we mean a graph $H$ which can
be obtained from $G$ via a sequence of edge contractions and vertex
deletions.
\end{defn}

This is different from the usual definition of a graph minor in
that we do not allow edge deletion, whose interaction with the
toric ideal is more complicated.

The rest of the paper is organized as follows.  In the next
section, we discuss computational results for graphs with few
vertices and suggest some conjectures based upon these data.  In
particular, we have computed the Markov width $\mu(G)$ for all
graphs on five vertices and many of the graphs on six vertices. In
the third section, we prove that the $n$-cycle and $K_{2,n}$ have
Markov width 4. Furthermore, we are able to explicitly describe
the moves needed in the Markov bases for these graph. The fourth
section is devoted to the inverse problem: that is, studying which
graphs may have Markov basis elements of a given degree.  We give
an algorithm which does not depend on computing $S$-pairs for
computing all the minimal generators of a given degree for the
ideals $I_G$.  As a consequence of our algorithm, we give a
combinatorial characterization of moves of degree two and three.
In the final section, we return to the study of forests and use
the reverse lexicographic Gr\"obner basis from Theorem
\ref{thm:sqfr2} to derive combinatorial formulae for the algebraic
degree of $I_G$ whenever $G$ is a forest.

\section{Graphs with Few Vertices}

In this section we discuss and display computational results about
the ideals of binary graph models.  In particular, we describe
generating sets for the ideals $I_G$ for all graphs $G$ with fewer
than five vertices and all the graphs on six vertices with at most
eight edges.  These computational results suggest many conjectures
and open problems which we describe at the end of the section.
All of our computations were carried out using the toric Gr\"obner
basis program \verb"4ti2" \cite{H} and the computational algebra
system \verb"Macaulay 2" \cite{M2}. We limit our description to
graphs which cannot be ``glued'' together from smaller graphs
based on the following definition and theorem.

\begin{defn}
Let $(V_1, S, V_2)$ be a partition of the vertex set of a graph
$G$ such that
\begin{enumerate}
\item there are no edges in $G$ between $V_1$ and $V_2$ and
\item $S$ is either the empty set or $S$ is a common vertex or
edge of the induced subgraphs $G_1$ and $G_2$ with vertex sets
$V_1 \cup S$ and $V_2 \cup S$, respectively.
\end{enumerate}

\noindent  Then $G$ is called {\bf reducible} with decomposition
$(V_1, S, V_2)$.
\end{defn}

\begin{thm}\label{reducible}\cite{DS2, HS}
Let $G$ be a graph which is reducible and $G_1$ and $G_2$ the
induced subgraphs arising from the vertex decomposition.  Then
there is a degree preserving operation with which one can build
generating sets and Gr\"obner bases for $I_G$ from generating sets
and Gr\"obner bases of $I_{G_1}$ and $I_{G_2}$. In particular,
$\mu(G) = \,{\rm max}(\mu(G_1),\mu(G_2))$.
\end{thm}

There are precisely one graph on three vertices, two graphs on
four vertices, six graphs on five vertices, and six graphs with
six vertices and at most eight edges which are not reducible.  We
will briefly describe these graphs and their Markov bases.

\subsection{Three and Four Vertices}

The only graph on three vertices which is not reducible is the
complete graph $K_3$. The Markov basis of $I_{K_3}$ consists of
the single degree four binomial
$$p_{000}p_{011}p_{101}p_{110} - p_{001}p_{010}p_{100}p_{111}.$$

The two irreducible graphs on four vertices are $C_4$ and $K_4$.
The graph $C_4$ is the four cycle with $E(C_4) = \{12,23,34,14\}$.
The Markov basis of $C_4$ consists of eight quadrics such as
$$ p_{0000}p_{0101} - p_{0001}p_{0100} $$
and eight quartics.  The complete graph $K_4$ has Markov basis
consisting of 20 moves of degree four and 40 sextic binomials such
as
$$ p_{0000}^2p_{0111}p_{1011}p_{1101}p_{1110} -
p_{0001}p_{0010}p_{0100}p_{1000}p_{1111}^2.$$

\subsection{Five Vertices}

There are six graphs on five vertices which cannot be decomposed
into subgraphs. These are the graphs we denote $C_5$, $K_{2,3}$,
$\widetilde{K_4}$, $SP$, $BP$, and $K_5$.  The graph $C_5$ is the
five-cycle $E(C_5) = \{12,23,34, 45,15\}$.  Its Markov basis
consists of 80 quadrics and 40 quartics.  The graph $K_{2,3}$ is
the complete bipartite graph $E(K_{2,3}) = \{13,14,15,23,
24,25\}$. Its Markov basis consists of 44 quadrics and 420
quartics.  The graph $\widetilde{K_4}$ is the graph obtained from
$K_4$ by subdividing an edge, $E(\widetilde{K_4}) =
\{12,15,23,24,34,35,45\}$.  The Markov basis for
$I_{\widetilde{K_4}}$ consists of 32 quadrics, 473 quartics, and
160 sextics.  The graph $SP$ is the edge graph of the square
pyramid, $E(SP) = \{12,13,15,23,24,34,35,45\}$.  The Markov basis
of $SP$ consists of 16 quadrics, 671 quartics, and 320 sextics.
The graph $BP$ is the edge graph of the bipyramid over a triangle,
$E(BP)  =  \{12,13,15,23,  \\ 24,25,34,35,45\}$.  Its Markov basis
consists of 8 quadrics, 436 quartics, and 2872 sextic binomials.
Finally, $K_5$ is the complete graph on five vertices.  The Markov
basis of $K_5$ consists of 260 degree four moves, 3952 sextic
binomials, 846 binomials of degree eight such as
$$
p_{00000}^3p_{01111}p_{10111}p_{11011}p_{11101}p_{11110} -
p_{00001}p_{00010}p_{00100}p_{01000}p_{10000}p_{11111}^3 $$

\noindent and 480 degree ten binomials like

$$
p_{00000}^2 p_{01111}^2 p_{10001}^2 p_{10010}^2p_{10100}p_{11000}-
p_{00010}p_{00101}p_{01001}p_{01110}p_{10000}^4p_{10011}p_{11111}.$$

\subsection{Six Vertices}

There are a total of 29 graphs on six vertices which are not
reducible.  We were able to compute Markov bases for the six
irreducible graphs on six vertices which have at most eight edges.
It remains a major computational challenge to determine Markov
bases of the other 23 irreducible graphs on six vertices.  The six
irreducible graphs on six vertices with less than nine edges will
be denoted $C_6$, $K_{2,4}$, $G_{129}$, $G_{151}$, $G_{153}$, and
$G_{154}$.

The graph $C_6$ is the six cycle with edge set $E(C_6) =
\{12,23,34,45,56,16\}$.  The graph $K_{2,4}$ is the complete
bipartite graph with edge set $E(K_{2,4}) =
\{13,14,15,16,23,24,25,26\}$.  The remaining graphs do not have
special names:  the labels we have chosen come from \cite{RW}.
These four graphs have edge sets $E(G_{129}) =
\{12,15,23,26,34,45,56\}$, $E(G_{151}) = \{12,14,23,26,34,36,45,46
\}$, $E(G_{153}) = \{12,15,16,23,24,45,46,56\}$, and $E(G_{154}) =
\{12,14,23,25,34,36,45,56\}$.  The data regarding the Markov bases
of these graphs as well as all the irreducible graphs on five and
fewer vertices is summarized in the following table. The columns
are labeled by the particular irreducible graph, the rows are
labelled by degree of minimal generators and the table entries are
the number of minimal generators of a given degree.

$$
\begin{array}{r||r|r|r|r|r|r|r|r|r|r|r|r|r|r|r|}
 & K_3 & C_4 & K_4 & C_5 & K_{2,3} & \widetilde{K_4} & SP & BP &
 K_5 & C_6 & K_{2,4} & G_{129} & G_{151} & G_{153} & G_{154} \\
 \hline 2 & 0 & 8 & 0 & 80 & 44 & 32 & 16 & 8 & 0 & 528 & 236 & 360& 280 & 320 & 256\\
 \hline 4 & 1 & 8 & 20 & 40 & 420 & 473 & 671 & 436 & 260 & 160 & 11696 & 2636 & 4949 & 4149 & 7784 \\
 \hline 6 & 0 & 0 & 40 & 0 & 0 & 160 & 320 & 2872 & 3952 & 0& 0 & 0 & 640 & 480 & 640 \\
 \hline 8 & 0 & 0 & 0 & 0 & 0 & 0 & 0 & 0 & 846 & 0& 0 & 0 & 0 & 0 & 0 \\
 \hline 10 & 0 & 0 & 0 & 0 & 0 & 0 & 0 & 0 & 480 & 0 & 0 & 0 & 0 & 0 & 0 \\
 \hline \hline \mbox{Total} & 1 & 16 & 60 & 120 & 464 & 665 & 1007 & 3316
 & 5538 & 688 & 11932 & 2996 & 5869 & 4949 &  8680\\
 \hline \mu(G) & 4 & 4 & 6 & 4 & 4 & 6 & 6 & 6 & 10 & 4 & 4 & 4 & 6 & 6  & 6 \\ \hline

\end{array}
$$

For these graphs, all generators are in even degree.  This is not
true in general, however, as we will demonstrate in Section 4.
Theorem~\ref{thm:sqfr2} characterizes graphs with $\mu(G)=2$ as
forests, but the next case is already quite interesting.

\begin{prob}
Characterize those graphs with Markov width $\mu(G) = 4$.
\end{prob}

In the next section we will show that cycles and the complete
bipartite graphs $K_{2,n}$ have Markov bases consisting of moves
of degree four or less, but from the data we see that this is not
yet a complete characterization.

A natural class of graphs which one would hope to
understand is planar graphs.  The data above suggest the following
optimistic conjecture.

\begin{conj}\label{planar}
There is a universal constant $C$ such that the Markov width
$\mu(G) \leq C$ whenever $G$ is a planar graph.  Even stronger, $C
= 6$.
\end{conj}

On the other hand, the data also suggest the following
conjecture.

\begin{conj}\label{treewidth}
The invariant $\mu(G)$ is a function only of the tree width of
$G$.
\end{conj}

The tree width is a topological invariant of a graph $G$ which is
equal to one less than the size of the largest clique in the
chordal graph containing $G$  which has the smallest maximal
clique.  For example, forests are precisely those graphs with tree
width zero or one, and indeed Theorem~\ref{thm:sqfr2} tells us
that these graphs all have Markov width two.
Conjecture~\ref{treewidth} also agrees with
Theorem~\ref{reducible}, since the tree width of a reducible graph
is the maximum of the tree widths of its components.

While the limited information we have suggests both
Conjecture~\ref{planar} and Conjecture~\ref{treewidth}, they
cannot both be true:  there are planar graphs with arbitrarily
large tree width.  For example, grid graphs can have arbitrarily
large clique size in their minimal chordal triangulations.  This
suggests another research problem.

\begin{prob}
Study the binary graph model $I_G$ for the family of $m \times n$
grid graphs.
\end{prob}

We do know that $\mu(G)$ can be arbitrarily large.  For example,
for the complete graph $K_{m}$ we can construct generators of
large degree.

\begin{prop}
The complete graph $K_m$ with $m \geq 3$ has Markov width
$\mu(K_m) \geq 2m -2$.
\end{prop}
\begin{proof}
It suffices to show that there is a minimal generator of $K_m$ of
degree $2m-2$.  For this consider the binomial

$$ p_{\mathbf{0}}^{m-2} \prod_{i=1}^m p_{\mathbf{1} - e_i} -
p_{\mathbf{1}}^{m-2} \prod_{i=1}^m p_{e_i}$$

\noindent where $\mathbf{0}$ is the string of all zeros,
$\mathbf{1}$ is the string of all ones, and $e_i$ is the $i$th
unit vector.  Then the monomials coming from the leading and
trailing terms are the only monomials which have the given image
under $\phi_{K_m}$.  Equivalently, the corresponding tables are
the only two tables which have these same fixed marginals under
$\pi_{K_m}$.  Since the leading and trailing terms have disjoint
support, this binomial must appear in every Markov basis for
$I_{K_m}$.
\end{proof}

Of course, this bound is already not tight for $m=5$,  where it
yields $\mu(K_5) \ge 8$ despite the fact that $K_5$ has Markov
width 10. In general we suspect that $\mu(K_m)$ grows
exponentially in $m$.

\section{Cycles and Bipartite Graphs}

In this section we confirm the observations from the second
section: the ideal of the cycle and the complete bipartite graph
$K_{2,n}$ are generated in degrees two and four.

\subsection{Cycles}

For ease of notation, we will represent a binomial such as
$p_{1011}p_{1110} - p_{1111} p_{1010}$ in {\bf tableau notation}
as

\[
\left[
\begin{array}{cccccc}
1 & 0 & 1 & 1 \\
1 & 1 & 1 & 0 \\
\end{array}
\right] - \left[
\begin{array}{cccccc}
1 & 1 & 1 & 1 \\
1 & 0 & 1 & 0 \\
\end{array}
\right].
\]

The tableau are obtained from a binomial by recording the indices
of each variable which appears in the monomial, repeating indices
when a variable appears to a power greater than one.  We say that
one binomial {\bf contains} another if it does so in the Graver
sense; that is, $p^u - p^v$ contains $p^a-p^b$ if $a \leq u $ and
$b \leq v$ componentwise.

We first prove the following theorem bounding the degree of
minimal generators for the $n$-cycle.

\begin{thm}\label{ncycle}
Let $C_n$ be the $n$-cycle graph. Then $\mu(C_n)=4$, and in particular $I_{C_n}$ is generated in degrees 2
and 4.
\end{thm}

\begin{proof}
We will start with an arbitrary binomial $f$ in the ideal, and
express it as a linear combination of elements either of lower
degree or of degree at most 4.

Given any binomial, take one variable from each monomial such that
the two variables chosen agree in first and last index. Our
strategy will be as follows: by adding multiples of ideal elements
of degree 4 and less, we will eventually obtain a binomial in
which both monomials have the same variable. Dividing out by this
variable (which clearly does not affect membership in $I_G$)
yields a binomial of lower degree which must still be in the
ideal, completing the proof.

We now start this process. We have a variable $p_{1?\cdots ?1}$ in
the first term of $f$, and a variable $p_{1?\cdots ?1}$ in the
second term. We wish to eliminate all disagreements between these
indices by ``moving'' the table entries corresponding to these
binary strings using binomials of degree 4 or less.

Consider any block of disagreements, in which, without loss of
generality, the indices of these variables look like
$(\cdots?10\cdots 01?\cdots)$ and $(\cdots?11\cdots11?\cdots)$. We
propose to add some multiple of an ideal element of degree $g$ at
most 4 so that the resulting binomial contains the two variables
in question, and is unchanged except that some of the
disagreements in the block have been removed. The two sets of
index strings in $g$ will agree on the portion of $C_n$ outside of
the block in question, counting the boundary elements.
Essentially, we are performing a local move by changing indices on
a subgraph of $C_n$. Let $\tilde{f}$ represent the image of $f$ in
this subgraph, i.e. under the map sending $p_{\cdots?I?\cdots}$ to
$p_I$, where $I$ is the index substring on the block we have.

Continuing in this manner, by induction on the number of
disagreements we eventually obtain a binomial for which the same
variable appears in both monomials, completing the proof. We now
construct the ideal element $g$ which we will add a multiple of.
We first construct the part $\tilde{g}$ which corresponds to the
block in question; in the tableaux that follow, we consider only
the indices corresponding to this block.

Because this element is in the ideal $I_{C_n}$, considering the
marginal in the first two directions, since an element in the
first term of $f$ has a $10$ marginal, so must an element in the
second term. So $\tilde{f}$ contains

\[
\left[
\begin{array}{cccccc}
1 & 0 & 0 & \cdots & 0 & 1 \\
\end{array}
\right] - \left[
\begin{array}{cccccc}
1 & 1 & 1 & \cdots & 1 & 1 \\
1 & 0 & ? & \cdots & ? & ? \\
\end{array}
\right] .
\]

Now, if any of the unspecified elements is 1, we let $\tilde{g}$
be the binomial which switches the intervening substrings, i.e.
something of the form

\[
\left[
\begin{array}{ccccc}
1 & 1 & 1 & 1 & 1 \\
1 & 0 & 0 & 1 & A \\
\end{array}
\right] - \left[
\begin{array}{ccccc}
1 & 0 & 0 & 1 & 1 \\
1 & 1 & 1 & 1 & A \\
\end{array}
\right],
\]
where $A$ is the remainder of the index string of the
element that the second term of $\tilde{f}$ contains starting with
$10$. We fill in the rest of both terms of $g$ as the two index
strings corresponding to the two variables contained in the second
term of $\tilde{f}$ are filled in, so that the two variables of
the first term of $g$ are contained in the second term of $f$.
Adding $g$ to $f$ then has the effect of eliminating some
disagreements between the two variables as desired, while leaving
everything unchanged outside the block in question.

Otherwise, all of the elements marked $?$ must be 0. We can apply
the same argument to the terminal string and to the other
binomials to obtain $g$ in all cases except where $\tilde{f}$
contains

\[
\left[
\begin{array}{ccccc}
1 & 0 & \cdots & 0 & 1 \\
1 & 1 & \cdots & 1 & 0 \\
0 & 1 & \cdots & 1 & 1 \\
\end{array}
\right] - \left[
\begin{array}{ccccc}
1 & 1 & \cdots & 1 & 1 \\
1 & 0 & \cdots & 0 & 0 \\
0 & 0 & \cdots & 0 & 1 \\
\end{array}
\right].
\]

In this case, because we have a $00$ marginal in the first two
coordinates of the second term of $\tilde{f}$, we must have one in
the first term. If that element contains any 1, by adding a
multiple of a binomial of degree 2 involving it and the third
element in the first term, and then another multiple of a binomial
of degree 2 involving the third element and the first element in
the first term, we can construct a $g$ essentially as before which
reduces disagreements. The only case where we cannot apply this
argument to this fourth element is when $\tilde{f}$ contains

\[
\left[
\begin{array}{ccccc}
1 & 0 & \cdots & 0 & 1 \\
1 & 1 & \cdots & 1 & 0 \\
0 & 1 & \cdots & 1 & 1 \\
0 & 0 & \cdots & 0 & 0 \\
\end{array}
\right] - \left[
\begin{array}{ccccc}
1 & 1 & \cdots & 1 & 1 \\
1 & 0 & \cdots & 0 & 0 \\
0 & 0 & \cdots & 0 & 1 \\
0 & 1 & \cdots & 1 & 0 \\
\end{array}
\right].
\]

In this case, we let $\tilde{g}$ be this binomial of degree 4,
corresponding to switching the middle substrings of all 0's and
all 1's. Again, we extend this to $g$ by copying the indices from
$f$ outside this block to the relevant variables, and we can add
this element to $f$ to eliminate this patch of disagreements
between the variables in question.  This completes the proof by
induction. Note that we did not use any elements of degree 3 in
this process, so $I_{C_n}$ is in fact generated in degrees 2 and
4.
\end{proof}

This theorem not only shows that the minimal generators are all of
degree 2 or 4, but it also gives a complete description of these
generators. The degree-2 generators come from separations of the
graph; we will prove a general statement characterizing degree-2
generators of graph ideals in Section 5.  As for the minimal
generators of degree four, we have the following categorization.

\begin{thm}
The minimal generators of degree 4 in the graph ideal $I_{C_n}$
are those elements of the form
\[
\left[
\begin{array}{cccc}
A_1 & 1 & B & 1 \\
A_2 & 1 & 1 -B & 0 \\
A_3 & 0 & 1 -B & 1 \\
A_4 & 0 & B & 0 \\
\end{array}
\right] - \left[
\begin{array}{cccc}
A_1 & 1 & 1-B & 1 \\
A_2 & 1 & B & 0 \\
A_3 & 0 & B & 1 \\
A_4 & 0 & 1-B & 0 \\
\end{array}
\right],
\]
where the columns correspond to $V_1, x_1, V_2$, and $x_2$, $V_1$
and $V_2$ are contiguous blocks of elements, and these elements in
this order comprise the $n$-cycle. Here, $1 -B$ represents the
opposite string of $B$.
\end{thm}

Note that these generators of degree 4 are very similar to the
generator of degree 4 in $K_3$, namely
\[
\left[
\begin{array}{ccc}
1 & 1 & 1 \\
1 & 0 & 0 \\
0 & 0 & 1 \\
0 & 1 & 0 \\
\end{array}
\right] - \left[
\begin{array}{ccc}
1 & 0 & 1 \\
1 & 1 & 0 \\
0 & 1 & 1 \\
0 & 0 & 0 \\
\end{array}
\right].
\]

We will prove a general similarity  theorem in this vein in
Section 4, when we classify generators of a given degree.

\subsection{Complete Bipartite Graphs}

Another nice class of models is the $K_{m,n}$ model, where
$K_{m,n}$ is the complete bipartite graph with partite sets of $m$
and $n$ vertices. We first prove the following theorem.

\begin{thm}\label{k2n}
The graph ideal for $G=K_{2,n}$ is generated in degrees 2 and 4
(for $n\ge 2$).
\end{thm}

\begin{proof}

As in the proof of Theorem~\ref{ncycle}, we use binomials of small
degree ($\le 4$) to transform a binomial in this ideal to one
whose two terms share a variable, completing the proof by
induction. Let the vertices of the two-element partite set be
$V=\{v_1,v_2\}$, and let the vertices of the $n$-element partite
set be $W=\{w_1,\ldots,w_n\}$.

For each monomial $M$, we define the submonomial $M_{ij}$, $i,j\in
\{0,1\}$, to be the product of the variables with $ij$ in the
index string corresponding to the two-element partite set; we will
write that index string first. For a monomial $M$, we define
$a_{ij}(M)$ to be the total degree of $M_{ij}$, and
$b_{ij,k,l}(M)$ to be the number of appearances of the digit $k$
in the $w_l$-position of the index strings of the variables in
$M_{ij}$. Here $k\in \{0,1\}$ and $l\in \{1,\ldots,n\}$.  In other
words, these function values enumerates the marginal in the
direction $(v_1,v_2,w_l)$ with the set values $(i,j,k)$.  We can
of course recover $a_{ij}(M) = b_{ij,0,l} + b_{ij,1,l}$ for any
$l$.

Then we have the following easy lemma.

\begin{lem}\label{deg2kmn}
If $M_1$ and $M_2$ are monomials such that $b_{ij,k,l}(M_1) =
b_{ij,k,l}(M_2)$ for all $i,j,k,l$, then their difference can be
expressed as a sum of multiples of quadratic elements of $I_G$.
\end{lem}

We do this simply by, for each $M_{ij}$, using quadratic
generators corresponding to the separation $(w_l, V, W\setminus
\{w_l\})$ to move around the $b_{ij,1,l}(M)$ 1's in the $l$th
column. Consequently, we need only to connect monomials with the
same marginals and different $F$-values. We introduce an
additional definition.

Given a monomial $M$, the function $c_{ij,l}(M)$ is defined to be
the subset of ${0,1}$ which appears in the $w_l$-position in the
variables of $M_{ij}$. Explicitly, this contains 0 when
$b_{ij,0,l}(M) > 0$, and 1 when $b_{ij,1,l}(M) > 0$.

We now unspool a series of moves designed to connect all of the
remaining $F$-values of monomials with the same marginals.

\begin{lem}\label{deg4sh1}
Suppose we have a monomial $M$ and a column $l$ such that $1\in
c_{01,l}(M), c_{10,l}(M)$ and $0\in c_{00,l}(M), c_{11,l}(M)$.
Then $M$ is equivalent by adding a multiple of a binomial of
degree 4 to a monomial with the following changes to the $a_I$'s
and $b_I$'s:
\[
\begin{array}{c}
+1: b_{01,0,l}, b_{10,0,l}, b_{11,1,l}, b_{00,1,l}, \\
-1: b_{01,1,l}, b_{10,1,l}, b_{11,0,l}, b_{00,0,l}. \\
\end{array}
\]
\end{lem}

\begin{proof}
This corresponds merely to adding a multiple of the degree 4
binomial

\[
\left[
\begin{array}{cccc}
1 & 1 & 1 & I_1\\
1 & 0 & 0 & I_2\\
0 & 0 & 1 & I_3\\
0 & 1 & 0 & I_4\\
\end{array}
\right] - \left[
\begin{array}{cccc}
1 & 0 & 1 & I_1\\
1 & 1 & 0 & I_2\\
0 & 1 & 1 & I_3\\
0 & 0 & 0 & I_4\\
\end{array}
\right],
\]
where the columns are $v_1, w_l$, and $v_2$, and all
other vertices in some order. This binomial comes from the minor
$K_3$ given by contracting all the $w_i$, $i\neq l$, into either
$v_1$ or $v_2$.
\end{proof}

\begin{lem}\label{deg2elt}

Suppose we have a monomial $M$ such that for each column $l$,
there exists some index $i_l$ such that $i_l\in c_{01,l}(M),
c_{10,l}(M)$. Then $M$ is equivalent by adding a multiple  of a
binomial of degree $2$ to a monomial with the following changes to
the $a_I$'s and $b_I$'s:
\[
\begin{array}{c}
+1: a_{11}, a_{00}, b_{11,i_l,l}, b_{00,i_l,l}, \\
-1: a_{10}, a_{01}, b_{01,i_l,l}, b_{10,i_l,l},
\end{array}
\]
for all $l\in \{1,\ldots,n\}$.

\end{lem}

\begin{proof}
If $I$ is the index string composed of the $i_l$, this corresponds
to adding a multiple of the degree 2 binomial

\[
\left[
\begin{array}{ccc}
1 & 1 & I\\
0 & 0 & I\\
\end{array}
\right] - \left[
\begin{array}{ccc}
1 & 0 & I\\
0 & 1 & I\\
\end{array}
\right],
\]
where the columns are indexed by $v_1$, $v_2$, and the $w_l$.
\end{proof}

Now, suppose we have two monomials with the same marginals, that
is the same image under $\phi_G$.  Add multiples of the binomials
from Lemma~\ref{deg4sh1} and Lemma~\ref{deg2elt} until one can no
longer apply these; since both increase $a_{11} + \sum b_{11,l}$,
one will not go around in circles. Our monomials $M$ and $N$ are
now in ``reduced'' form, in the sense that neither move can be
applied. We break the situation down into cases.

{\it Case 1. Suppose that $M_{11}$ and $N_{11}$ are both not equal
to $1$, so that both $M$ and $N$ have an entry which is $11$ in
the $(v_1,v_2)$ direction.}

If, for each $l$, there exists an index $i_l$ such that $i_l\in
c_{11,l}(M)$ and $i_l\in c_{11,l}(N)$, then, as desired, we can
simply extract the variable $p_{11,(i_l)}$ from both $M$ and $N$;
in other words, for both $M$ and $N$, we can find a monomial with
the same values of $a$ and $b$ containing this variable.

If this is not the case, then there exists an $l$ for which
without loss of generality $c_{11,l}(M) = \{1\}$ and $c_{11,l}(N)
= \{0\}$. Looking at the $(v_1,w_l)$ marginal, there exists at
least one marginal $11$ because of the first condition;
consequently, there must exist at least one of these marginals in
$N$. Since $c_{11,l}(N)=\{0\}$, the only other option is that
$1\in c_{10,l}(N)$. Similarly, considering the $(v_2,w_l)$
marginal, we must have $1\in c_{01,l}(N)$.

Now, we have $b_{11,0,l}(N) = k>0$. Since $b_{11,0,l}(M)=0$,
looking at the $10$-count in the $(v_1,w_l)$ direction, we must
have $b_{10,0,l}(M) = b_{10,0,l}(N)+k$. However, looking at the
$00$-count in the $(v_2,w_l)$ direction, it now follows that we
must have $b_{00,0,l}(N) = b_{00,0,l}(M)+k$, and in particular
$0\in c_{00,l}(N)$. This is a contradiction, since we can now
apply a move as in Lemma~\ref{deg4sh1} to $N$, contradicting the
assumption that $N$ is reduced.

{\it Case 2. Exactly one of $M_{11}$ and $N_{11}$ is equal to 1.}

Suppose without loss of generality that $M_{11}\neq 1$ and
$N_{11}=1$. Take any $i_l\in c_{11,l}(M)$; this $i_l$ must be in
$c_{10,l}(N)$ and $c_{01,l}(N)$. This means that we can apply a
move as in Lemma~\ref{deg2elt} to $N$, again contradicting the
hypothesis that $N$ is reduced.

{\it Case 3. Both $M_{11}$ and $N_{11}$ are empty.}

In this case, it follows immediately that
$b_{10,i,l}(M)=b_{10,i,l}(N)$ for all $i$ and $l$ by considering
the $(v_1,w_l)$ marginals equal to $(1,i)$. If $M_{10}\neq 1$,
this means that we can find an $i_l$ for all $l$ such that this
number is nonzero, and we can then pull the corresponding variable
out of both $M$ and $N$. Similarly, if $M_{01}\neq 1$, we can find
a shared variable there, and if both of these are 1 then applying
the same argument to $M_{00}=M$ and $N_{00}=N$ finishes the job.

Our litany of cases has come to an end, completing the proof of
Theorem~\ref{k2n}. Note again that we have not only  shown that
$\mu(K_{2,n})=4$, but also given an explicit generating set in
degrees 2 and 4 for $I_{K_{2,n}}$.
\end{proof}

For $K_{m,n}$ where $m,n>2$, the answer is less clear. The
statement and proof of Theorem~\ref{k2n} indicate that for $m$
fixed, as $n$ gets large, the maximum degree of an element in the
Markov basis of the graph ideal of $K_{m,n}$ stabilizes. The
degree, on the other hand, certainly goes up as $\text{min}(m,n)$
does; for instance, there is an element of degree $2m$ in the
Markov basis of $K_{m,m}$, and we have the following result.

\begin{prop}
Fix $m\ge 2$. Then for $n \geq {m \choose 2} 2^{m-2}$, there   is
an element of degree $2^{m-1}$ in the graph ideal of $K_{m,n}$.
\end{prop}

\begin{proof}
Let the vertices of $K_{m,n}$ be $\{v_1,\ldots,v_m\}\cup \{w_I\}$,
where $I=(i_1,\ldots,i_m)$ is an index string of length $m$,
consisting of precisely two 1's, and some number of 0's and 2's.
There are precisely ${m \choose 2} 2^{m-2}$ such strings.

Specify the marginals as follows: between $v_j$ and $w_I$, insist
upon $i_j$ marginals of $11$, $2-i_j$ marginals of $01$,
$2^{m-2}-i_j$ marginals of $10$, and $2^{m-2}-2+i_j$ marginals of
$00$. What this means is that exactly two variables with
coordinate $w_I$ equal to 1 occur, and that the sum of the
$v$-coordinates (considered as vectors) is precisely $I$; it
furthermore specifies that each of 0 and 1 occurs $2^{m-2}$ times
in each $v_j$-coordinate.

For each $w_I$, there are only two ways to express the vector  $I$
as the sum of two 0-1 vectors. Since we consider all index strings
$I$, we obtain that for each diamond in the natural Boolean
partial order of binary strings of length $m$, either the top and
bottom elements are in the set of $v$-coordinates of table
entries, or the middle two entries are. By an easy induction, it
follows that the set of $v$-coordinates, which numbers only
$2^{m-1}$, must consist of either all strings with an even number
of 1's, or all strings with an odd number of 1's. From here, we
can easily compute the $w$-coordinates of each of these entries.

These resulting tables are the only two which satisfy  these
marginals, and thus their difference, an element of degree
$2^{m-1}$, must be in the Markov basis of $I_{K_{m,n}}$ as
desired.   For all $n \geq {m \choose 2} 2^{m-2}$ there is a move
of degree $2^{m-1}$ by Corollary \ref{cor:minor}.
\end{proof}

We suspect that the following conjecture, an  extension of the
result for $K_{2,n}$, holds.

\begin{conj}
The graph ideal for $G=K_{m,n}$ is generated in degree at most
$2^{\min\{m,n\}}$.
\end{conj}

\section{Combinatorial Classification of Minimal Generators of Low Degree}

In this section, we give algorithms for computing  all generators
of a given degree in the graph ideal $I_G$. For degrees two and
three, we give an explicit combinatorial characterization of these
generators, giving a generating set which generates $I_G$ in
degree less than or equal to 3; for arbitrary degree $d$, we
categorize these generators as pullbacks of a distinguished
generator in the graph ideal of a fundamental graph $X_d$. The key
lemma is the following, relating generators in $I_G$ to generators
in a minor of $G$.

\begin{lem}\label{gphom}
Let $G$ be a graph, and let $f=\Pi p_{I_j} - \Pi p_{I_k}$ be a
binomial contained in $I_G$. Then we have the following.

(a) If $v_i$ corresponds to a column where all the index strings
$I_j$ have the same value, then $f$ is a minimal generator if and
only if $\tilde{f}$ is a minimal generator of the graph ideal
$G\backslash v_i$, where $\tilde{f}$ is the natural image of $f$
with the column $v_i$ deleted from each index string.

(b) If $v_i$ and $v_j$ are adjacent and correspond to columns
where for each index string $I_j$ or $I_k$, the value of that
string in each column is identical, then $f$ is a minimal
generator if and only if $\tilde{f}$ is a minimal generator, where
$\tilde{f}$ is the natural image of $f$ with the two columns $v_i$
and $v_j$ fused. Here, $\tilde{f}$ is an element of the graph
ideal of $G$ with those two vertices fused.

(c) Suppose $v_i$ and $v_j$ are any two vertices and correspond to
columns where for each index string $I_j$ or $I_k$, the value of
that string in each column is identical. In this case, if $f$ is a
minimal generator of the graph ideal of $G$, then $\tilde{f}$ is a
minimal generator of the graph ideal of $G$ with those two
vertices fused.
\end{lem}

\begin{proof}
In each case, a decomposition of $\tilde{f}$ into generators of
lower degree can be lifted via the obvious method to a
decomposition of $f$. In (a), this is simply inserting the shared
value of $v_i$ into each index string to form a valid index string
for $G$; in (b) and (c), this is simply duplicating the value of
each index string in the obvious manner.

In cases (a) and (b), any decomposition of $f$ must necessarily
satisfy the property that each binomial used has the property in
question, by considering in (a) any marginal containing $v_i$ and
in (b) the marginal corresponding to the edge $v_iv_j$. Therefore,
a decomposition of $f$ naturally yields a decomposition of
$\tilde{f}$, so if $\tilde{f}$ is a minimal generator $f$ must be
also.
\end{proof}

This lemma has a corollary legitimizing our notion of minor.

\begin{cor} \label{cor:minor}
If $H$ is a minor of $G$, then $\mu(H)\le \mu(G)$.
\end{cor}

\begin{proof}
It suffices to show $\mu(H)\le \mu(G)$ if $H$ is obtained  from
$G$ by a single vertex deletion or edge contraction. However, if
it is obtained by a vertex deletion, then by part (a) of
Lemma~\ref{gphom}, every minimal generator of a given degree in
$I_H$ lifts to a minimal generator of the same degree in $I_G$.
Similarly, part (b) of Lemma~\ref{gphom} guarantees that
$\mu(H)\le \mu(G)$ if $H$ is obtained from $G$ via an edge
contraction.
\end{proof}

A natural extension of this is the following, which  agrees with
the data in Section 2, but which we have been unable to prove.

\begin{conj}
If $H$ is obtained from $G$ by deleting an edge, then $\mu(H)\le \mu(G)$.
\end{conj}

This set of minimal generators comes with a group action. In
particular, the group $(\ZZ/2)^n$ acts naturally on
$\mathbb{C}[p_I]$, via the element $(c_1,\ldots,c_n)$ sending a
variable $p_{i_1\cdots i_n}$ to $p_{j_1,\ldots,j_n}$, where $j_r =
i_r + c_r$; the sum is evaluated in $\ZZ/2$. This action consists
merely of flipping 0's and 1's in some positions.

Furthermore, the automorphism group
$\text{Aut}(G)$ acts naturally on $\mathbb{C}[p_I]$ as well, by
permuting the indices according to the permutation of the vertices
of the graph, so we have a natural action of $(\ZZ/2)^n \oplus
\text{Aut}(G)$ on $\mathbb{C}[p_I]$. This action maps $I_G$ onto
itself; we make the following natural definition.

\begin{defn}
Two generators are {\bf equivalent} if they lie in the same orbit
of $\mathbb{C}[p_I]$ under the action of $(\ZZ/2)^n \rtimes
\text{Aut}(G)$.
\end{defn}

If two generators of graph ideals $I_{G_1}$ and $I_{G_2}$ reduce
to equivalent generators in a basic graph $H$ by means of the
above manipulations, we say that they are {\bf weakly similar}. If
they furthermore reduce to equivalent generators using only
manipulations of type (a) and (b), we say that they are {\bf
strongly similar}. We are now prepared to define  the object
pivotal in our categorization of generators of degree $d$.

\begin{defn}\label{fungraph}
Fix a degree $d\ge 2$. The {\bf fundamental graph}  $X_d$ has
vertex set $(S_i,T_i)$, where $S_i$ and $T_i$ are subsets of
$\{1,\ldots,d\}$ with cardinalities $|S_i| = |T_i| \leq d/2$, and
if $|S_i|\,=d/2$ then $1\in S_i$. Two vertices $(S_1,T_1)$ and
$(S_2,T_2)$ are connected by an edge if $|S_1\cap
S_2|\,=\,|T_1\cap T_2|$.
\end{defn}

To this fundamental graph is associated a distinguished  element
of $I_{X_d}$.

\begin{defn}
Fix a degree $d\ge 2$. Then the {\bf distinguished generator}
$f_d\in I_{X_d}$ is the binomial

\[
\left[
\begin{array}{c}
I_1\\
\cdots\\
I_d\\
\end{array}
\right] - \left[
\begin{array}{c}
J_1\\
\cdots\\
J_d\\
\end{array}
\right],
\]

where $I_{ji}=1$ if $j\in S_i$ and 0 otherwise,  and similarly
$J_{ji}=1$ if $j\in T_i$ and 0 otherwise.
\end{defn}

It is clear that the distinguished generator is  actually in
$I_{X_d}$, since for all adjacent $v_i$ and $v_j$, the number of
11-marginals in $I$ is equal to $S_i\cap S_j$, while the number of
11-marginals in $J$ is equal to $T_i\cap T_j$. By definition of
$X_d$, these are equal, and furthermore, the number of
11-marginals determines the numbers of all other marginals (along
with the numbers of 1's in each column of $I$ and $J$, which are
of course identical.)

\begin{exmp}
The fundamental graph $X_2$ has two vertices $(1,1)$ and $(1,2)$
which are not connected by an edge.  The distinguished generator
of $I_{X_2}$ is the binomial

$$ \left[ \begin{array}{cc} 1 & 1 \\ 0 & 0 \end{array} \right] -
\left[ \begin{array}{cc} 1 & 0 \\ 0 & 1  \end{array}\right]$$

\noindent which is just as $2 \times 2 $ determinant.  The
fundamental graph $X_3$ has nine vertices  which are $(1,1),
(1,2), \ldots , (3,3)$.  Two vertices $(i,j)$ and $(k,l)$ are
connected if and only if $i \neq k$ and $j \neq l$.  Each vertex
has degree four:  $X_3$ is the edge graph of the 4-polytope
$\Delta_2 \times \Delta_2$, the product to two triangles pictured
in Figure 1. Note the six triangular prisms which appears as
minors of $X_3$. The fundamental generator of $I_{X_3}$ is the
binomial

$$ \left[ \begin{array}{ccccccccc}
1 & 1 & 1 & 0 & 0 & 0 & 0 & 0 & 0 \\
0 & 0 & 0 & 1 & 1 & 1 & 0 & 0 & 0 \\
0 & 0 & 0 & 0 & 0 & 0 & 1 & 1 & 1 \end{array} \right] - \left[
\begin{array}{ccccccccc}
1 & 0 & 0 & 1 & 0 & 0 & 1 & 0 & 0 \\
0 & 1 & 0 & 0 & 1 & 0 & 0 & 1 & 0 \\
0 & 0 & 1 & 0 & 0 & 1 & 0 & 0 & 1 \end{array} \right].$$
\end{exmp}

\begin{figure}
\begin{center}\includegraphics{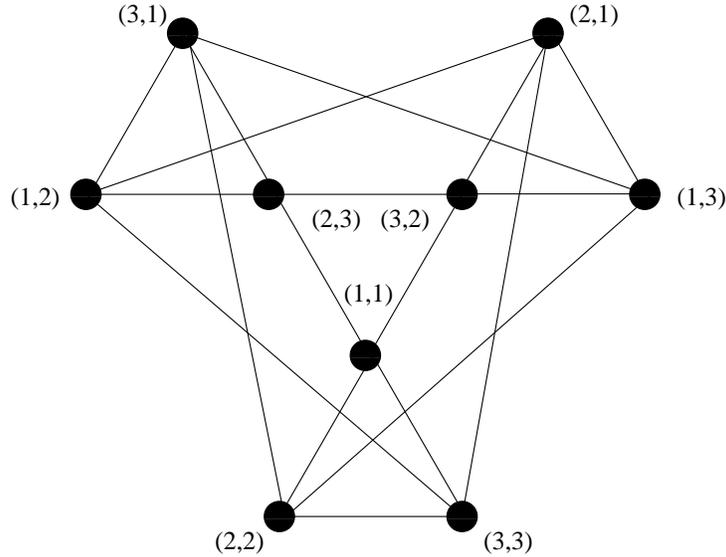} \caption{The fundamental graph
$X_3$}
\end{center}
\end{figure}

With these definitions, we can formulate the main  theorem of this
section, categorizing all generators of degree $d$.

\begin{thm}\label{degreed}
Let $G$ be any graph. Then the minimal generators of  degree $d$
in $I_G$ can be enumerated by the following procedure.
\begin{enumerate}
\item Consider all graph homomorphisms $\phi$ from minors $H$ of $G$ to
$X_d$.
\item Determine if the natural image of the fundamental generator $\widetilde{f_d}$
is a minimal generator on the image subgraph.
\item If so, consider the pullback $\phi^\star \widetilde{f_d}$, a binomial of degree
$d$ in $I_G$.  We can pull this back to $H$  using
Lemma~\ref{gphom} part (b) and (c), and then to $G$ using
Lemma~\ref{gphom} part (a) and (b).
\end{enumerate}
\end{thm}

A graph homomorphism $G\rightarrow H$ is simply a map $\phi$  from
vertices of $G$ to vertices of $H$ such that if $u$ and $v$ are
adjacent in $G$, $\phi(u)$ and $\phi(v)$ are adjacent in $H$. Note
that since we have used Lemma \ref{gphom} part (c), the pullback
will not always be a minimal generator but the set of moves
calculated in this way contains all minimal generators $I_G$ of
degree $d$.

\begin{proof}
Suppose we have a generator $f$ of degree $d$ in $I_G$,  written
in tableau notation as

\[
\left[
\begin{array}{c}
I_1\\
\cdots\\
I_d\\
\end{array}
\right] - \left[
\begin{array}{c}
J_1\\
\cdots\\
J_d\\
\end{array}
\right],
\]
where $I$ and $J$ are 0-1 matrices. Each column of $I$ has the
same number of 1's as the corresponding column of $J$. If this
number is either 0 or $d$, we can delete that vertex via
Lemma~\ref{gphom} (a) to obtain an equivalent generator in a minor
of $G$, proving the theorem by induction.

Therefore, by flipping 0 with 1 if necessary, we can assume that
each column has at most $d/2$ 1's, and if it has $d/2$ 1's then
$I_i$ has a 1 in that column. To each vertex $v_i$ associate the
pair $(S_i,T_i)$, where $j\in S_i$ if $I_j$ has a 1 in the column
corresponding to $v_i$, and $j\in T_i$ if $J_j$ has a 1 in that
column. If two adjacent vertices have the same pair, we can
contract the edge between them via Lemma~\ref{gphom} (b) to again
obtain an equivalent generator in a minor of $G$.

If this is not the case, then we claim that the map given by
sending $v_i$ to the vertex $(S_i,T_i)$ in $X_d$ is a graph
homomorphism. Indeed, all that we need to check is that if $v_i$
and $v_j$ are connected by an edge, $S_i\cap S_j = T_i\cap T_j$
(and it is not the case that $(S_i,T_i)=(S_j,T_j)$, which is true
since no two adjacent vertices have the same pair.) But this must
be true, since each is just the number of 11-marginals along the
edge $v_iv_j$ in the corresponding tables $I$ and $J$, which are
equal since $f$ is in $I_G$.

Furthermore, this map fuses two vertices $v_i$ and $v_j$ only if
the corresponding columns of $I$ and of $J$ are identical.
Therefore, by Lemma~\ref{gphom} (c), the image of $f$ in $I_R$,
where $R$ is the image subgraph of $X_d$, must be a minimal
generator of $X_{I_R}$. However, this generator is precisely the
natural image of the distinguished generator in the graph ideal
$I_R$. If this is irreducible, then $f$ is produced by the
procedure in Theorem~\ref{degreed}, which we have just done in
reverse. If not, then $f$ cannot be a minimal generator by the
contrapositive of Lemma~\ref{gphom}.
\end{proof}

In this manner, we have reduced the computation of all generators
of degree $d$ to the process of determining which images of the
distinguished generator in subgraphs of $X_d$ are minimal
generators, and of enumerating graph homomorphisms from $G$ to
$X_d$. While the problem of computing graph homomorphisms is a
difficult one, we can use symmetry techniques to greatly aid us in
many cases. Theorem~\ref{degreed} also divulges which generators
are weakly similar: those whose corresponding graph homomorphisms
have images which are isomorphic subgraphs of $X_d$.

We now apply Theorem~\ref{degreed} to degrees 2 and 3. Consider
first generators of degree 2. The fundamental graph $X_2$ has two
vertices, $(1,1)$ and $(1,2)$, which are not connected by an edge.
Given a minor of $G$, it has a homomorphism onto $X_2$ if and only
if there are no edges between the vertices mapped to $(1,1)$ and
the vertices mapped to $(1,2)$. In this case the homomorphisms
correspond to partitions of the vertices into $V_1$ and $V_2$.
Extracting the definition of minor yields the following Corollary.

\begin{cor}
Let $G$ be a graph.  Then equivalence classes of generators of
degree 2  of $I_G$ correspond to partitions $V_1\cup V_2\cup V_3$,
where there are no edges between $V_1$ and $V_2$. The generator
corresponding to this is precisely

\[
\left[
\begin{array}{cc}
1 0 1 \\
0 1 1\\
\end{array}
\right] - \left[
\begin{array}{cc}
1 1 1\\
0 0 1\\
\end{array}
\right],
\]
where the three columns correspond to $V_1, V_2$, and $V_3$.
\end{cor}

Readers familiar with the study of graphical models and their
induced independence statements will recognize this as a theorem
which says that the only independence statements induced by a
graphical model are the global independence statements.  See, for
example \cite{L}.

Next, we turn our attention to the case of generators of degree 3,
using the methods of Theorem~\ref{degreed} to obtain a
combinatorial classification of all such cubic minimal generators
of $I_G$. The fundamental graph $X_3$ has nine vertices
$\{(1,1),\ldots,(3,3)\}$, with $(i,j)$ connected to $(k,l)$ if
$i\neq k$ and $j\neq l$.  By direct computation, the image of the
fundamental generator $d_3$ in $I_R$, $R$ a subgraph of $X_3$, is
a minimal generator if and only if $R$ contains one of the six
triangular prism subgraphs of $X_3$. Therefore, we can classify
cubic generators by finding all homomorphisms from minors of $G$
to the fundamental graph $X_3$ whose image contains a triangular
prism. In particular, if no such homomorphism exists, then $I_G$
cannot have a generator of degree 3. This yields the following
corollary.

\begin{cor}\label{no3kn}
For all $n$, the graph ideal $I_{K_n}$ has no minimal generators of degree 3.
\end{cor}
\begin{proof}
The only minors of $K_n$ are copies of $K_m$ for $m\le n$.
Furthermore, all homomorphic images of complete graphs are again
complete graphs (indeed, of the same degree.) The only complete
graphs occurring in $X_3$ are $K_2$ and $K_3$, neither of which
contains a triangular prism.
\end{proof}

Corollary~\ref{no3kn} shows why edge deletion  does not behave
well with respect to the graph ideal $I_G$ and its Markov basis;
this operation obviously can introduce elements of new degrees,
since there exist graphs with cubic minimal generators, and every
graph can be obtained from a complete graph by edge deletion. The
same technique can be used to show that there are no minimal
generators of degree 5 in the graph ideal $I_{K_n}$, since the
largest clique in $X_5$ has size 5, and the distinguished
generator cannot be minimal in the graph ideal of these subgraphs
since in fact there are no minimal generators of degree 5 in $K_i$
for $i\le 5$.

In addition to the description of  Theorem~\ref{degreed}, we can
obtain a more straightforward combinatorial characterization of
minimal generators of degree 3 in graph ideals. We start with a
pair of definitions.

\begin{defn}
Let $G$ be any graph. Then the {\bf 3-coloring graph} $C_3(G)$ has
vertices equal to the set of proper 3-colorings of $G$. Two
3-colorings are connected by an edge if one can be obtained from
the other by the following (reversible) procedure: pick a color
$i\in \{1,2,3\}$, pick a connected component of the induced
subgraph consisting of all vertices with colors not equal to $i$,
and switch the other two colors on this component.
\end{defn}

\begin{defn}
A graph $G$ is {\bf 3-rigid} if the coloring space $C_3(G)$ is disconnected.
\end{defn}

The crucial example of a 3-rigid graph is the same triangular
prism that arose in the analysis of degree 3 via
Theorem~\ref{degreed}. This graph is 3-rigid, since there are only
two proper 3-colorings up to permutation of the colors: one
triangle $(v_1,v_2,v_3)$ is colored $(1,2,3)$, and the other
triangle $(w_1,w_2,w_3)$ must be colored either $(2,3,1)$ or
$(3,1,2)$. Here, the vertex labels are chosen so that $v_i$ is
adjacent to $w_i$. It is easy to check that these colorings lie in
different connected components of $C_3(G)$; the connected
component of each consists of the colorings obtained from it by
permutation of the colors.

We now present a complete  description of cubic generators of $G$
based on these combinatorial objects.

\begin{prop}\label{3rigid}
Let $G$ be any graph. Then  we can enumerate the (equivalence
classes of) cubic minimal generators of the graph ideal $I_G$ as
follows.

(a) Find all the 3-rigid minors of $G$.

(b) For each such 3-rigid minor $H$, consider all of the connected
components of the 3-coloring graph of $H_3$. For each connected
component, pick a representative 3-coloring $R_i$.

(c) For each $R_i$, $i>1$, take the element
\[
f= \left[
\begin{array}{ccc}
a_1\\
a_2\\
a_3\\
\end{array}
\right] - \left[
\begin{array}{ccc}
b_1\\
b_2\\
b_3\\
\end{array}
\right],
\]
where $a_{jk}=1$ if either the vertex $k$ of $G$ is deleted in
obtaining $H$ or its image in $H$ is colored with color $j$ in
$R_1$, and $b_{jk}=1$ if either $k$ is deleted or its image in $H$
is colored with color $j$ in $R_i$.

Each of these elements $f$ is a minimal generator of $I_G$ in
degree 3, and these elements together with the quadratic elements
described above generate $I_G$ up to degree 3.
\end{prop}

\begin{proof}
Suppose we have a cubic minimal generator of $I_G$, given by

\[
f= \left[
\begin{array}{ccc}
a_1\\
a_2\\
a_3\\
\end{array}
\right] - \left[
\begin{array}{ccc}
b_1\\
b_2\\
b_3\\
\end{array}
\right].
\]

We construct two 3-colorings $C_A$ and $C_B$ of a minor of $G$ as
follows. For the vertex $v_j$, consider the $j$-th column of $A$
and $B$. If this consists of all 1's or all 0's, delete $v_j$. If
not, exchange 0's and 1's if necessary so that it has exactly one
1. Then, in $C_A$ (resp. $C_B$), color $v_j$ with the position in
which this 1 appears in $A$ (resp. B).

Next, if two adjacent vertices have the same color in $A$
(equivalent to having the same color in $B$ by counting
11-marginals along this edge), then contract that edge. What
remains is two proper 3-colorings of a minor $H$, and when these
3-colorings are converted to a binomial in $I_G$ via the method in
part (c), we recover precisely the element $f$. By
Lemma~\ref{gphom}, the image $\tilde{f}$ in $I_H$ is a minimal
generator if and only if $f$ is a minimal generator of $I_G$.

We claim that the image $\tilde{f}$ is a minimal  generator of
$I_H$, i.e. inexpressible as a sum of multiples of quadrics, if
and only if $C_A$ and $C_B$ are in different connected components
of $C_3(H)$. Indeed, consider a multiple of a binomial generator,
say
\[
g= \left[
\begin{array}{ccc}
a_1\\
a_2\\
a_3\\
\end{array}
\right] - \left[
\begin{array}{ccc}
c_1\\
c_2\\
a_3\\
\end{array}
\right].
\]

If we convert the two monomials of $g$ into colorings, the same
set of vertices will have color 3 in these colorings. Therefore,
the difference between these colorings consists of changing colors
from 1 to 2 or vice versa. The only way this can be done while
preserving the properness of the coloring is if the change
consists of switching the colors 1 and 2 on some connected
components of the induced subgraph of $H$ consisting of vertices
not colored 3. Therefore, when we add a multiple of a binomial
generator, we end
 1
up with an element corresponding to a 3-coloring in the same
component of $C_3(H)$, and any two elements connected by an edge
in $C_3(H)$ differ by a multiple of a binomial generator.

Consequently, the colorings $C_A$ and $C_B$ are in  the same
component of $C_3(H)$ if and only if $\tilde{f}$ is not
expressible as the sum of multiples of binomials in $I_H$, which
is equivalent to $\tilde{f}$ being a minimal generator of $I_H$,
as desired.
\end{proof}

Classifying 3-rigid graphs is an interesting problem;  the graph
$C_3(G)$ has been studied in connection with the problem of
picking a random 3-coloring of a graph~\cite{Vigoda}. Indeed, the
flip interchanging two colors on a connected component is
precisely the move used in the Wang-Swendsen-Koteck\'{y} algorithm
to pick a random $k$-coloring of a graph, and this scheme has ties
to mathematical physics~\cite{WSK}. There are simple operations to
produce 3-rigid graphs from other ones, but the triangular prism
seems to be the only 3-rigid graph without a proper 3-rigid minor.

This method, paralleling Theorem~\ref{degreed}, can be  extended
to higher degrees. However, it rapidly becomes unwieldy, as the
vertices can now be colored with sets of $d/2$ colors, and the
moves are more complicated, consisting of all moves keeping one of
the colors fixed.  For instance, the fundamental graph $X_4$ has
34 vertices and understanding the homomorphisms to this graph
seems difficult.

\section{Algebraic Degree of Forests}

A recent series of results gave a very thorough description of the
family of ideals of decomposable models \cite{D, GMS, T}. A
special case of these results is the following fundamental
theorem.

\begin{thm}\label{thm:sqfr}
The ideal $I_G$ is minimally generated by quadrics if and only if
$G$ is a forest.  In this case, the set of quadratic squarefree
binomials in $I_G$ forms a Gr\"obner basis with respect to the
reverse lexicographic term order with $p_{0 \cdots 0} \prec p_{0
\cdots 0 1} \prec \ldots \prec p_{1 \cdots 1}$.
\end{thm}

Sturmfels \cite{St2} posed the natural follow-up problem of
calculating the degree of the toric ideal $I_G$.  The degree of
the toric ideal is interesting in statistics because it provides a
natural upper bound for the maximum likelihood degree of the toric
ideal \cite{St2}.  In this section we give combinatorial formulae
for the degree of the graph model ideal $I_G$ when $G$ is a
forest.  As the maximum likelihood degree of a forest is always 1,
we see that the degree can be arbitrarily far from the maximum
likelihood degree.  To perform these degree computations, we first
recall a result about the degree of a general toric ideal. This
result can be found in \cite{St}.

\begin{thm}
Let $A$ be a $d \times n$ matrix whose toric ideal $I_A$ is
homogeneous in the standard $\ZZ$-grading.  Then the degree of the
ideal $I_A$ ($=$ degree of the projective toric variety defined by
$I_A$) is equal to the normalized volume of the lattice polytope
$Q = \mathrm{conv}(A)$.
\end{thm}

Henceforth, the normalized volume of a lattice polytope $Q$ will
be denoted $Vol(Q)$.  This theorem reduces the problem of
calculating degree to computing the normalized volume of polytope.
We now record some some basic facts about the polytope $P_G$ when
the underlying graph is a forest.

\begin{lem}
The polytope $P_G$ has dimension equal to the sum of the number of
vertices and the number of edges of the graph (this is true for
any graph).  There are precisely $4 \cdot | E(G) |  + 2 \cdot |
Iso(G) | $ facets of $P_G$ when $G$ is a forest.  If the variables
in marginal space are labelled $y^{(j,k)}_{i_j, i_k}$ for the
variables coming from an edge and $y^{(l)}_{i_l}$ for the
variables coming from an isolated vertex, then the facets are
given by the inequalities

$$ y^{(j,k)}_{i_j, i_k} \geq 0 \mbox{  and   } y^{(l)}_{i_l} \geq 0.$$

\end{lem}

\begin{proof}
The dimension formula appears in \cite{HS}.  The polyhedral
results are a direct consequence of the closed form expressions
for maximum likelihood estimates in decomposable models in
\cite{L}.
\end{proof}

The following lemma implies that to compute the degree of $I_G$
when $G$ is a forest, one need only describe a formula which is
valid for trees.  Furthermore, this lemma is important for
carrying out the recursive computations implied by Theorem
\ref{thm:vol}.

\begin{lem}\label{lem:prod}
Let $G$ be a graph with a partition of the vertices $\{V_1,V_2\}$
such that there is no edge in $G$ between the $V_1$ and $V_2$. Let
$G_1$ and $G_2$ be the corresponding induced subgraphs.  Let $d_1$
and $d_2$ be the corresponding dimensions of the polytopes
$P_{G_1}$ and $P_{G_2}$; that is, $d_i = |V_i| + |E(G_i)|$ .  Then
we have

$$\label{eqn:prod} Vol(P_G) = { d_1 + d_2 \choose d_1} \cdot Vol(P_{G_1}) \cdot
Vol(P_{G_2}). $$

\end{lem}

\begin{proof}
With these restrictions on the graph $G$, we have $P_G = P_{G_1}
\times P_{G_2}$.  Equation \ref{eqn:prod} is then the usual
formula for the normalized volume of the direct product in terms
of the normalized volumes of the components of the product.
\end{proof}

We now come to the main theorem of this section.

\begin{thm}\label{thm:vol}
Let $G$ be a tree with $n$-vertices.  Then the degree of the toric
ideal $I_G$ can be calculated by the formula

$$deg(I_G) = \frac{1}{2} \sum_{e \in E(G)} deg(I_{G \setminus e})
$$

\noindent where the notation $G \setminus e$ denotes the graph $G$
with the edge $e$ removed.
\end{thm}

\begin{proof}
As previously indicated, we prove the theorem by calculating the
volume of the corresponding polytope $P_G$.  Theorem
\ref{thm:sqfr} implies that the pulling triangulation of $P_G$
induced by the reverse lexicographic term order above is
unimodular. This in turn, implies that the normalized volume of
$P_G$ is equal to the sum of the normalized volumes of the facets
of $P_G$ which do not contain the vertex indexed by the variable
$p_{1\cdots 1}$ (see \cite{St} for all the definitions and
relevant theory). This is where the polyhedral description of
$P_G$ when $G$ is a forest becomes essential.  We see from the
polyhedral characterization that there are exactly $n-1$ facets of
$P_G$ which do not contain this ``last'' vertex, and that they are
indexed by the edges of $G$.  We will show that the normalized
volume of the facet $F_e$ of $P_G$ which is indexed by the edge
$e$ has volume precisely $\frac{1}{2} Vol(P_{G \setminus e})$
which will complete the proof.

There are two cases to consider:  either the edge in question has
one node a leaf or not (the case of the graph on two vertices with
a lone edge is clear, by a direct calculation). We will handle the
two cases separately.

\emph{Case 1}:  We may suppose without loss of generality that our
edge $e$ is $\{1,2\}$, the vertex 1 is a leaf and the vertex $2$
has the edge $\{2,3\}$ incident to it.  Then the matrix whose
columns correspond to the vertices of $P_G$ has the top eight rows
which look like

$$ \left( \begin{array}{cccccccc}
1 & 1 & 0 & 0 & 0 & 0 & 0 & 0 \\
0 & 0 & 1 & 1 & 0 & 0 & 0 & 0 \\
0 & 0 & 0 & 0 & 1 & 1 & 0 & 0 \\
0 & 0 & 0 & 0 & 0 & 0 & 1 & 1 \\
1 & 0 & 0 & 0 & 1 & 0 & 0 & 0 \\
0 & 1 & 0 & 0 & 0 & 1 & 0 & 0 \\
0 & 0 & 1 & 0 & 0 & 0 & 1 & 0 \\
0 & 0 & 0 & 1 & 0 & 0 & 0 & 1 \end{array} \right) $$

\noindent with this $8 \times 8$ block repeated $2^{n-3}$ times
across the first eight rows.  We claim that the facet which
corresponds to the inequality $y^{(1,2)}_{1,1} \geq 0$
(corresponding to the fourth row of the above matrix) has volume
equal to $\frac{1}{2} Vol(P_{G \setminus \{1,2\}})$.  Note the the
vertices of $P_G$ which lie on this facet are precisely the $3
\cdot 2^{n-2}$ vertices of $P_G$ which have a zero in the fourth
row.  First we show that this facet is naturally isomorphic to a
sub-configuration of $P_{G \setminus \{1,2\}}$. Consider the
matrix $A_{G \setminus \{1,2\}}$ whose columns give the vertices
of $P_{G \setminus \{1,2\}}$.  This matrix has two fewer rows than
the matrix $A_G$ above and is almost the same: its first six rows
look like

$$ \left( \begin{array}{cccccccc}
1 & 1 & 1 & 1 & 0 & 0 & 0 & 0 \\
0 & 0 & 0 & 0 & 1 & 1 & 1 & 1 \\
1 & 0 & 0 & 0 & 1 & 0 & 0 & 0 \\
0 & 1 & 0 & 0 & 0 & 1 & 0 & 0 \\
0 & 0 & 1 & 0 & 0 & 0 & 1 & 0 \\
0 & 0 & 0 & 1 & 0 & 0 & 0 & 1
 \end{array} \right).$$

\noindent To show the natural isomorphism mentioned above, it
suffices to show that there is a unimodular linear transformation
from the first six columns of the first matrix above to the first
six columns of the second matrix above.  Such a linear
transformation is obtained by replacing the first row by the sum
of the first and second rows, and then deleting the second and
fourth rows.  We can delete the second and fourth rows because
they are linear combinations of other rows and hence do not change
the polyhedral description.  Such a transformation is clearly
unimodular.

Now that we have shown that our configuration of $3 \cdot 2^{n-2}$
points sits naturally inside $P_{G \setminus \{1,2\}}$, we wish to
compute its volume.  For this, we show that there is a hyperplane
which divides $P_{G \setminus \{1,2\}}$ into two congruent pieces,
one of which is the convex hull of our new configuration of $3
\cdot 2^{n-2}$ points.  This hyperplane is given by the equation
$$ y^{(1)}_0 - y^{(1)}_1 + y^{(2,3)}_{0,0} + y^{(2,3)}_{0,1} - y^{(2,3)}_{1,0}
-y^{(2,3)}_{1,1} = 0.$$

\noindent  Note that exactly $2^{n-1}$ vertices of $P_{G \setminus
\{1,2\}}$ lie on this hyperplane (these are the ones corresponding
to the middle four columns of the submatrix  of $A_{G \setminus
\{1,2\}}$ displayed above) and the remaining $2^{n-1}$ vertices
are split equally on each side of the hyperplane.  In particular,
all of the vertices from our configuration of $3 \cdot 2^{n-2}$
points lie on the nonnegative side of this hyperplane and the
remaining $2^{n-2}$ points are on the negative side.  Furthermore,
there is a natural reflexive symmetry across this hyperplane.  To
complete the proof, we must show more: not only is the point
configuration naturally ``cut in half'' by this hyperplane, but so
is the polytope $P_{G \setminus \{1,2\}}$. This follows from a
direct computation with the eight points listed above.  We
performed the computation using the program \verb"PORTA" \cite{C}.

\emph{Case 2}: The argument is very similar to case 1;   we will
sketch the relevant details.  We may assume that our edge is the
edge $\{2,3\}$.  Since neither 2 nor 3 is a leaf we may assume $G$
also contains the edges $\{1,2\}$ and $\{3,4\}$.  With these
conditions, the first 12 rows of our matrix looks like

$$ \left( \begin{array}{cccccccccccccccc}
1 & 1 & 1 & 1 & 0 & 0 & 0 & 0 & 0 & 0 & 0 & 0 & 0 & 0 & 0 & 0 \\
0 & 0 & 0 & 0 & 1 & 1 & 1 & 1 & 0 & 0 & 0 & 0 & 0 & 0 & 0 & 0 \\
0 & 0 & 0 & 0 & 0 & 0 & 0 & 0 & 1 & 1 & 1 & 1 & 0 & 0 & 0 & 0 \\
0 & 0 & 0 & 0 & 0 & 0 & 0 & 0 & 0 & 0 & 0 & 0 & 1 & 1 & 1 & 1 \\
1 & 1 & 0 & 0 & 0 & 0 & 0 & 0 & 1 & 1 & 0 & 0 & 0 & 0 & 0 & 0 \\
0 & 0 & 1 & 1 & 0 & 0 & 0 & 0 & 0 & 0 & 1 & 1 & 0 & 0 & 0 & 0 \\
0 & 0 & 0 & 0 & 1 & 1 & 0 & 0 & 0 & 0 & 0 & 0 & 1 & 1 & 0 & 0 \\
0 & 0 & 0 & 0 & 0 & 0 & 1 & 1 & 0 & 0 & 0 & 0 & 0 & 0 & 1 & 1 \\
1 & 0 & 0 & 0 & 1 & 0 & 0 & 0 & 1 & 0 & 0 & 0 & 1 & 0 & 0 & 0 \\
0 & 1 & 0 & 0 & 0 & 1 & 0 & 0 & 0 & 1 & 0 & 0 & 0 & 1 & 0 & 0 \\
0 & 0 & 1 & 0 & 0 & 0 & 1 & 0 & 0 & 0 & 1 & 0 & 0 & 0 & 1 & 0 \\
0 & 0 & 0 & 1 & 0 & 0 & 0 & 1 & 0 & 0 & 0 & 1 & 0 & 0 & 0 & 1 \\
\end{array} \right) $$

\noindent  with this block repeated $2^{n-4}$ times.  We wish to
show that the facet defined by the inequality $y^{(2,3)}_{1,1}
\geq 0 $ (corresponding to the eighth row in the above matrix) has
volume equal to $\frac{1}{2} Vol(P_{G \setminus \{2,3\}})$.  The
vertices which lie on this facet are precisely the $3 \cdot
2^{n-2}$ vertices with a zero in the eighth row in this matrix
representation.  First we show that this facet naturally appears
as a sub-configuration of $P_{G \setminus \{2,3\}}$.  This can be
achieved by applying a unimodular transformation to the
configuration:  the key point is that once we restrict attention
to the facet, the middle 4 rows of the above configuration can be
written as linear combinations of the other rows and hence are
redundant in terms of the polyhedral description.

Now we show that there is a hyperplane which divides the polytope
$P_{G \setminus (2,3)}$ in half.  This is just the hyperplane
given by
$$y^{(1,2)}_{0,0} - y^{(1,2)}_{0,1} + y^{(1,2)}_{1,0} - y^{(1,2)}_{1,1} +
y^{(3,4)}_{0,0} + y^{(3,4)}_{0,1} - y^{(3,4)}_{1,0} -
y^{(3,4)}_{1,1} = 0 .$$

\noindent  Note that our configuration of $3 \cdot 2^{n-2}$ points
are precisely the points on the nonnegative side of this
hyperplane.  Furthermore, there is a natural reflexive symmetry
across the hyperplane.  A direct calculation shows that this
hyperplane not only separates the point configuration, but also
divides the polytope into two symmetric pieces with the the
desired integral points as vertices.  Thus we deduce the desired
equation of volumes.
\end{proof}

For some special families of trees we use this recurrence relation
to deduce simple formulae for the degree.  These appear in the
following corollaries.

\begin{cor}
Let $K_{1,n}$ denote a star graph with $n$ leaves.  Then
$deg(I_{K_{1,n}}) = (n!)^2$.
\end{cor}

\begin{proof}
Removing any edge of the graph $K_{1,n}$ produces the graph
consisting of the disjoint union of a $K_{1, n-1}$ and an isolated
point.  Hence from theorem \ref{thm:vol} and lemma \ref{lem:prod}
we deduce that

$$deg(I_{K_{1,n}}) = n \cdot \frac{1}{2} {2n \choose 1} \cdot
deg(I_{K_{1,n-1}}) = n^2 \cdot deg(I_{K_{1,n-1}}).$$

\noindent  Since $deg(I_{K_{1,1}}) =1$ we have the desired result.
\end{proof}

\begin{cor}
Let $T_n$ denote the graph of the $n$-chain and $d_n =
deg(I_{T_n})$.  Then $d_n$ satisfies the recurrence

$$ \label{eqn:rec} d_{n+1} = \frac{1}{2} \sum_{i=1}^n  {2n \choose 2i -1} d_i
d_{n+1-i}$$

\noindent with $d_1 = 1$.  Furthermore, we have the equality of
generating functions

$$\sum_{n=1}^{\infty} \frac{d_n}{(2n-1)!} x^{2n-1} = \sqrt{2}
\tan(\frac{x}{\sqrt{2}}).$$
\end{cor}

\begin{proof}
The recurrence relation \ref{eqn:rec} follows immediately from the
formula in Theorem \ref{thm:vol} and by applying Lemma
\ref{lem:prod} to the graph consisting of two disjoint chains of
length $i$ and $n+1-i$.

To deduce the equality of generating functions, let $y =
\frac{1}{2} \sum_{n=1}^{\infty} \frac{d_n}{(2n-1)!} x^{2n-1}$. The
recurrence relation implies that $2y' -2 = y^2$.  Solving the
differential equation yields the desired formula.
\end{proof}

The recurrence relation and generating function in the case of the
$n$-chain also appears in a paper by Poupard \cite{Po}, but we do
not know how to show that the objects we are counting (simplices
in a regular unimodular triangulation) are in bijection with the
objects she was counting (types of binary trees).

\section{Acknowledgements}

Mike Develin was supported in part by an NSF Graduate Research
Fellowship and by the American Institute of Mathematics.  Seth
Sullivant was supported by an NSF Graduate Research Fellowship.


\begin{thebibliography}{99}

\bibitem{C} T.~Christof.
\verb"PORTA": a polyhedron representation transformation
algorithm, version 1.3.1 available at
\verb"ftp://elib.zib.de/pub/mathprog/polyth/index.html".

\bibitem{DS} P.~Diaconis and B.~Sturmfels.
Algebraic algorithms for sampling from conditional distributions.
\emph{Annals of Statistics}, 26 (1998), 363--397.

\bibitem{D} A.~Dobra.
Markov bases for decomposable graphical models.  \emph{Bernoulli},
9 (2003), no. 6, 1--16.

\bibitem{DS2} A.~Dobra and S.~Sullivant.
A divide-and-conquer algorithm for generating Markov bases of
multi-way tables.  To appear in \emph{Computational Statistics}.

\bibitem{Dyer} M.~Dyer, L.A.~Goldberg, C.~Greenhill, M.~Jerrum, M.~Mitzenmacher. An extension of path
coupling and its application to the Glauber dynamics for graph colorings. {\em SIAM J. Comput.},
30 (2001), 1962--1975.

\bibitem{GMS} D.~Geiger, C.~Meek, and B.~Sturmfels.
On the toric algebra of graphical models. Manuscript, 2002.

\bibitem{M2}
D. Grayson and M. Stillman. \verb"Macaulay 2": a software system
for research in algebraic geometry.  Available at
\verb"http://www.math.uiuc.edu/Macaulay2/".

\bibitem{H} R.~Hemmecke.
\verb"4ti2": computation of Hilbert bases, Graver bases, toric
Gr\"{o}bner bases and more.  Available at
\verb"http://www.4ti2.de/".

\bibitem{HS}
S. Ho\c{s}ten and S. Sullivant. Gr\"{o}bner bases and polyhedral
geometry of cyclic and reducible models,  \emph{Journal of
Combinatorial Theory: Series A}, 100 (2002), no. 2, 277--301.

\bibitem{L} S.~Lauritzen. {\em Graphical Models}. Oxford
University Press, New York, 1996.

\bibitem{Po}
C. Poupard. Deux proprietes des arbres binaires ordonnes stricts,
\emph{European J. Combin.}, 10 (1989), 369--374.

\bibitem{RW}
R. Read and R. Wilson. {\em An Atlas of Graphs}.  Oxford
University Press, New York, 1999.

\bibitem{St}
B.~Sturmfels. \emph{Gr{\"o}bner Bases and Convex Polytopes},
American
  Mathematical Society, Providence, 1995.

\bibitem{St2}
B.~Sturmfels.
\newblock {\em Solving Systems of Polynomial Equations}, volume~97 of {\em CBMS
  Regional Conference Series in Mathematics}.
\newblock American Mathematical Society, Providence, 2002.


\bibitem{T} A.~Takken.
{\em Monte Carlo Goodness-of-Fit Tests for Discrete Data}. PhD
thesis,  Stanford University,  1999.

\bibitem{Vigoda} E. Vigoda. Improved bounds for sampling colorings,
\emph{J. Mathematical Physics}, 41 (2000) no. 3, 1555--1569.

\bibitem{WSK} J.S. Wang, R.H. Swendsen, and R. Koteck\'{y}.
Antiferromagnetic Potts models. \emph{Phys. Rev. Lett}, 63 (1989), no. 2, 109--112.

\end{thebibliography}
\end{document}